\documentclass[reqno,11pt]{amsart}

\IfFileExists{mymtpro2.sty}{%
  \usepackage[subscriptcorrection]{mymtpro2}
}{}

\usepackage{a4,amssymb}
\usepackage{mathrsfs}
\usepackage{graphicx}

\newtheorem{theorem}{Theorem}

\newtheorem{thm}{Theorem}[section]
\newtheorem{lemma}{Lemma}[section]
\newtheorem{follow}{Corollary}[section]

\theoremstyle{definition}

\newcommand{\bel}{\begin{equation} \label}
\newcommand{\ee}{\end{equation}}
\newcommand{\R}{\mathbb{R}}

\newcommand{\C}{\mathbb{C}}

\theoremstyle{remark}
\newtheorem{remark}[theorem]{Remark}
\newtheorem{myremarks}[theorem]{Remarks}

    {\end{nummer}\end{myremarks}}

\newcounter{numcount}
\newcommand{\labelnummer}{\mbox{\normalfont (\roman{numcount})}}%

\makeatletter

\newenvironment{nummer}%
  {\let\curlabelspeicher\@currentlabel%
    \begin{list}{\labelnummer}%
      {\usecounter{numcount}\leftmargin0pt%
        \topsep0.5ex\partopsep2ex\parsep0pt\itemsep0ex\@plus1\p@%
        \labelwidth2.5em\itemindent3.5em\labelsep1em%
      }%
    \let\saveitem\item%
    \def\item{\saveitem%
      \def\@currentlabel{{\upshape\curlabelspeicher}$\,$\labelnummer}}%
    \let\savelabel\label%
    \def\label##1{\savelabel{##1}%
      \@bsphack%
        \ifmmode\else%
          \protected@write\@auxout{}%
          {\string\newlabel{##1item}{{\labelnummer}{\thepage}}}%
        \fi%
      \@esphack%
    }%
  }{\end{list}}%

\renewcommand{\appendix}{\def\thesection{\textsc{Appendix}}}

 \let\leq\le
 \let\geq\ge

\DeclareMathOperator{\tr}{tr\kern1pt}

\newcommand{\supp}{\mathop\mathrm{supp}\nolimits}

\makeatletter

\newif\ifper\pertrue
\def\per{.}

\def\bti{\@ifnextchar[\bbti\bbbti}
\def\bbti[#1]#2{#2, #1.}
\def\bbbti#1{#1.}

\def\z{\@ifnextchar[\zz\zzz}
\def\zz[#1]#2#3#4#5{\perfalse\emph{#2} \textbf{#3}, #4 (#5) [#1]}
\def\zzz#1#2#3#4{\emph{#1} \textbf{#2}, #3 (#4)\ifper\per\fi\pertrue}

\def\pub{\@ifstar\pubstar\pubnostar}
\def\pubnostar{\@ifnextchar[\@@pubnostar\@pubnostar}
\def\@@pubnostar[#1]#2#3#4{#2, #3, #4, #1\ifper\per\fi\pertrue}
\def\@pubnostar#1#2#3{#1, #2, #3\ifper\per\fi\pertrue}
\def\pubstar[#1]#2#3#4{\perfalse #2, #3, #4 [#1]\pertrue}

\makeatother

\newcommand{\beq}{\begin{equation}}
\newcommand{\eeq}{\end{equation}}
\newcommand{\ba}{\begin{array}}
\newcommand{\ea}{\end{array}}
\newcommand{\bea}{\begin{eqnarray}}
\newcommand{\eea}{\end{eqnarray}}
\newcommand{\beas}{\begin{eqnarray*}}
\newcommand{\eeas}{\end{eqnarray*}}

{

\def\P{I\kern-.30em{P}}
\def\E{I\kern-.30em{E}}
\renewcommand{\E}{\mathbb{E}\mkern2mu}
\renewcommand{\P}{\mathbb{P}}

\begin{document}

\title[Multiplicity of Hodge eigenvalues]{The multiplicity of eigenvalues of the Hodge Laplacian on $5$-dimensional compact manifolds}

\author[M.\ E.\ Gier]{Megan E.\ Gier}
\address{
    Department of Mathematics,
    Grove City College,
     Grove City, PA 16127 }
\email{MEGier@gcc.edu}

\author[P.\ D.\ Hislop]{Peter D.\ Hislop}
\address{Department of Mathematics,
    University of Kentucky,
    Lexington, Kentucky  40506-0027}
\email{peter.hislop@uky.edu}

\thanks{Both authors were partially supported by NSF grant DMS 11-03104 during the time this work was done. This paper is partly based on the dissertation submitted by the first author in partial fulfillment of the requirements for a PhD at the University of Kentucky.}

\begin{abstract}
We study multiplicity of the eigenvalues of the Hodge Laplacian on smooth, compact Riemannian manifolds of dimension five for generic families of metrics. We prove that generically
the Hodge Laplacian, restricted to the subspace of co-exact two-forms, has nonzero eigenvalues of multiplicity two.
The proof is based on the fact that Hodge Laplacian restricted to the subspace of co-exact two-forms is minus the square of the Beltrami operator, a first-order operator. We prove that for generic metrics the spectrum of the Beltrami operator is simple. Because the Beltrami operator in this setting is a skew-adjoint operator, this implies the main result for the Hodge Laplacian.
\end{abstract}

\maketitle \thispagestyle{empty}


\tableofcontents

\vspace{.2in}



\section{Statement of the problem and results}\label{sec:introduction}
\setcounter{equation}{0}

The multiplicity of the $L^2$-eigenvalues of the Laplacian $\Delta_g \geq 0$ on a smooth compact manifold $(M,g)$ is linked with the symmetry of the manifold. Generally speaking, the multiplicity of an eigenvalue is reduced under perturbations of the Laplacian.
In the seventies, Uhlenbeck \cite{uhlenbeck1} and Albert \cite{albert1} studied this question for generic classes of metric and potential perturbations. For a Riemannian manifold $(M, g_0)$, Uhlenbeck proved that a generic, local perturbation of the metric $g_0 \rightarrow g_0 + \delta g$, with support $\delta g \subset U \subset M$, an open set, removes all multiplicities. That is, the eigenvalues of $\Delta_{g_0 + \delta g}$ are  simple (that is, have multiplicity one) for a generic set of perturbations $\delta g$ supported in $U \subset M$.

In light of Uhlenbeck's result for the Laplace operator on functions, one might wonder if the nonzero eigenvalues of the Hodge Laplacian $\Delta_g^{(k)}$ acting on $k$-forms might likewise be simple for a residual set of metrics.  Soon after Uhlenbeck published her theorem, Millman \cite{Millman} noted that on a manifold of even dimension $2n$, the McKean-Singer t$\acute{\mbox{e}}$lescopage theorem \cite{bgm1} implies that all the nonzero eigenvalues of the Hodge Laplacian acting on $n$-forms (forms of middle dimension) have even multiplicity.  While Millman's observation precludes a general extension of Uhlenbeck's theorem to the Hodge Laplacian, it is possible for analogues to hold under appropriate hypotheses.

In 2012, Enciso and Peralta-Salas \cite{EPS} proved that on a closed 3-manifold, there exists a residual set of $C^r$ metrics, $r \geq 2$, such that the nonzero eigenvalues of the Hodge Laplacian $\Delta_g^{(k)}$, for $0\leq k\leq 3$, all have multiplicity 1.  They structure their proof around the study of the Beltrami operator $*_gd$ restricted to co-exact 1-forms, which they show to have simple spectrum by a similar transversality theory argument as employed by Uhlenbeck. The Beltrami operator $*_gd$ restricted to co-exact 1-forms is self-adjoint and its square, on the same subspace, is the Hodge Laplacian $\Delta_g^{(1)}$, restricted to this invariant subspace. Consequently, the Hodge Laplacian
restricted to this subspace also has simple nonzero eigenvalues. This fact, when combined with the Hodge decomposition and Uhlenbeck's theorem for the Laplace operator acting on 0-forms (functions), and Hodge duality, allow Enciso and Peralta-Salas to conclude their simplicity result for $\Delta_g^{(1)}$.  The generic simplicity of the nonzero spectrum of the Hodge Laplacian acting on $k$-forms for $0\leq k\leq 3$ follows from Uhlenbeck's theorem for $k=0$, their result for $k=1$, and Hodge duality for $k=2$ and $k=3$.

In this paper, we extend the method centered on the Beltrami operator, as introduced by Enciso and Peralta-Salas \cite{EPS}, to study the generic nonzero eigenvalue multiplicities of the Hodge Laplacian on closed 5-manifolds.  In particular, we will prove that for a residual set of $C^r$ metrics, for any $r \geq 2$, the nonzero eigenvalues of the Hodge Laplacian $\Delta_g^{(2)}$ acting on co-exact 2-forms have multiplicity 2. Instead of transversality, we employ the direct perturbation theory method used by Albert \cite{albert1} (also used by Colin de Verdi\`ere \cite{cdv1}).

In order to state the main theorem, we recall the de Rham complex of real differential forms over $M$.
The de Rham complex for $(M,g)$ consists of the spaces $\Lambda^k(M)$ of smooth $k$-forms on $M$ and the differential maps $d:\Lambda^k(M)\to\Lambda^{k+1}(M)$ for $k=0,\ldots,n$. Each $\Lambda^k(M)$ is a pre-Hilbert space with inner product given by
\begin{eqnarray}\label{eq:ip}
(u,v)_g &=& \int_M u \wedge {(*_g v)} \hspace{.1in}\mbox{for }u,v\in\Lambda^k(M),
\end{eqnarray}
where $\wedge$ is the wedge product and $*_g:\Lambda^k(M)\to\Lambda^{n-k}(M)$ is the Hodge star operator. We denote the closure of $\Lambda^k(M)$ in the related norm by $L^2 (M, \Lambda^k)$. In the discussion of the Beltrami operator in section \ref{sec:beltrami-ev1}
we will work with complex-valued forms that we denote by $\Lambda_{\C}^k(M)$. In this case, the form ${(*_g v)}$
in the inner product \eqref{eq:ip} is replaced by its complex conjugate denoted by $\overline{(*_g v)}$. The adjoint of $d$ with respect to this inner product is the codifferential operator $\delta_g:\Lambda^{k+1}(M)\to\Lambda^k(M)$.  Our primary operator of interest is the Hodge Laplacian, the second order differential operator given by $\Delta_g^{(k)}=d\delta_g+\delta_gd$, acting on its natural domain in $L^2 (M, \Lambda^k)$.

The operators $\Delta_g^{(k)}$, $d$, and $\delta_g$ allow us to define the following subspaces of $\Lambda^k(M)$.  The space of harmonic $k$-forms on $M$ is
$$\mathcal{H}^k(M)=\{u\in\Lambda^k(M)|\,\Delta_g^{(k)}u=0\},$$
the space of exact $k$-forms is
$$d\Lambda^{k-1}(M)=\{u\in\Lambda^k(M)|\,u=dv\mbox{ for some }v\in\Lambda^{k-1}(M)\},$$
and the space of co-exact $k$-forms is
$$\delta_g\Lambda^{k+1}(M)=\{u\in\Lambda^k(M)|\,u=\delta_gw\mbox{ for some }w\in\Lambda^{k+1}(M)\}.$$
The Hodge Decomposition Theorem guarantees that any $k$-form can be uniquely written as the sum of a harmonic form, an exact form, and a co-exact form:

\begin{thm}\label{thm:Hodge}\cite{Morita}
On an oriented compact Riemannian manifold $(M,g)$, the space $\Lambda^k(M)$ can be decomposed as
$$\Lambda^k(M)=\mathcal{H}^k(M)\oplus d\Lambda^{k-1}(M)\oplus\delta_g\Lambda^{k+1}(M).$$ The space of harmonic forms $\mathcal{H}^k(M)$
is finite dimensional.
\end{thm}

The result extends to an orthogonal decomposition of $L^2 (M, \Lambda^k)$. If $H^1 (M, \Lambda^k)$ is the Sobolev space of $k$-forms,
then  $L^2 ( M, \Lambda^k ) = \mathcal{H}^k(M)\oplus d H^1(M,\Lambda^{k-1})\oplus \delta_g H^1(M,\Lambda^{k+1})$, see, for example, \cite[Theorem 1.5.2]{gilkey}.

The Beltrami operator $*_gd$ maps $k$-forms to $n-k-1$-forms, with the ranks of the forms coinciding precisely when $n=2k+1$.  In particular, the manifold must be of odd dimension. In the case studied by Enciso and Peralta-Salas with $n=3$, the Beltrami operator maps $1$-forms to $1$-forms. The spectrum of the Hodge Laplacian restricted to exact $1$-forms follows from Uhlenbeck's analysis of the spectrum of the Laplace-Beltrami operator on $0$-forms since the exact $1$-forms have the form $d f$. On co-exact $1$-forms, the Hodge Laplacian equals a phase factor times the square of the Beltrami operator. Hence, by the Hodge decomposition, the spectrum of the Hodge Laplacian on $1$-forms is determined by the Beltrami operator. By Hodge duality, this determined the spectrum of the Hodge Laplacian on $2$-forms.

The next dimension for which the Beltrami operator may be used to study the spectrum of the Hodge Laplacian is $n = 5$. In this case, the Beltrami operator
maps $2$-forms to $2$-forms.
In particular, the square of the Beltrami operator acting on co-exact $2$-forms is minus the Hodge Laplacian acting on co-exact $2$-forms. Consequently, the Beltrami operator may be used to study the spectrum of the Hodge Laplacian restricted to the invariant subspace of co-exact $2$-forms.

\begin{thm}\label{thm:hodge-main1}
 Let $M$ be a closed, 5-dimensional Riemannian manifold. Let $r$ be an integer with $r\geq 2$.  There exists a residual subset $\Gamma$ of the space of all $C^r$ metrics on $M$ such that, for all $g\in\Gamma$, the nonzero eigenvalues of the Hodge Laplacian $\Delta_g^{(2)}$ acting on co-exact 2-forms have multiplicity 2.
\end{thm}

Our proof of Theorem \ref{thm:hodge-main1} centers on an investigation of the Beltrami operator $*_gd$.  Using perturbation theory inspired by Albert \cite{albert1}, and a density argument of Colin de Verdi\`ere \cite{cdv1}, we will show that for a residual set of metrics, the Beltrami operator restricted to co-exact 2-forms has only simple eigenvalues.  We will then explore the relationship between the spectrum of the Beltrami operator, a skew-adjoint operator, and that of the Hodge Laplacian on co-exact 2-forms. In particular, the origin of the generic multiplicity two of eigenvalues is the skew-adjointness of the Beltrami operator on 2-forms. This means the eigenvalues of the Beltrami operator are pure imaginary and the real and imaginary parts of the complex eigenforms give rise to independent real eigenforms of the Hodge Laplacian.   The main result follows from this.

\subsection{The meaning of generic}

In this article, the terms \emph{generic} and \emph{generic property} mean the following.
Let $X$ be a topological space. A set $\mathcal{G} \subset X$ will be called \emph{residual} or \emph{generic} in $X$ if it is a dense $G_\delta$-set. That is, $\mathcal{G} = \cap_{j=1}^\infty G_j$, where each $G_j \subset X$ is dense and open in $X$.  A property that is true for a residual subset of a topological space $X$ is called \emph{generic}.

\subsection{Discussion of the Beltrami and Hodge operators}

The Beltrami operator may be used to study the eigenvalues of the Hodge Laplacian restricted to co-exact $k$-forms
only for certain pairs $(n,k)$ of dimension $n$ of the manifold and rank $k$ of the forms.
Before narrowing our focus to co-exact 2-forms on a 5-manifold, we consider the more general properties of the Beltrami operator acting on $k$-forms on an $n$-dimensional manifold. Since the Beltrami operator is the composition of $*_g$ and $d$, the operator is an isomorphism between $\delta_g\Lambda^{k+1}(M)$ and $\delta_g\Lambda^{n-k}(M)$, that is, the spaces of real co-exact $k$-forms and co-exact $(n-k-1)$-forms. The Beltrami operator may be extended to complex-valued forms by linearity. The extended Beltrami operator $*_gd:\delta_g\Lambda_{\mathbb{C}}^{k+1}(M)\to\delta_g\Lambda_{\mathbb{C}}^{n-k}(M)$ is also an isomorphism.

\begin{lemma}\label{lem:coexact}
Let $M$ be an $n$-manifold.  Then $$\Delta_g^{(k)}=(-1)^{nk+1}(*_gd)^2$$ when restricted to co-exact, real or complex, $k$-forms.
\end{lemma}

\textit{Proof.} If $\omega\in\delta_g\Lambda_{\mathbb{C}}^{k+1}(M)$, then
	$\Delta_g^{(k)}\omega = \delta_gd\omega$. In terms of the Hodge star operator, the co-differential operator $\delta_g$ is
 $ \delta_g = (-1)^{n(k+1)+1}*_g d *_g$. Using this, we find
 $$
 \Delta_g^{(k)}\omega = (-1)^{n(k+2)+1}(*_gd*_g)d\omega = (-1)^{nk+1}(*_gd)^2\omega.
$$
The same calculation holds on $\delta_g\Lambda^{k+1}(M)$. \hfill$\Box$

Lemma \ref{lem:coexact} implies that when restricted to co-exact forms, the Hodge Laplacian is given by $\Delta_g^{(k)}=(*_gd)^2$ if $n$ and $k$ are both odd; otherwise $\Delta_g^{(k)}=-(*_gd)^2$.  The parity of $n$ and $k$ also determine whether the Beltrami operator is self-adjoint or skew-adjoint.

\begin{lemma}\label{lem:adjoint}
Let $M$ be an $n$-dimensional manifold, $\omega\in H^1(M,\Lambda_{\mathbb{C}}^k)$, and $\eta\in H^1(M,\Lambda_{\mathbb{C}}^{n-k-1})$.  Then $$( *_gd\omega, \eta)_g=(-1)^{nk+1}( \omega,*_gd\eta)_g.$$
\end{lemma}

This result indicates that the Beltrami operator is self-adjoint if $(n,k)$ are both odd and skew-adjoint otherwise. Combining this with the mapping properties of the Beltrami operator, we make the following conjecture concerning
the generic multiplicities of the nonzero eigenvalues of the Hodge Laplacian on odd dimensional manifolds:
The nonzero eigenvalues of the Hodge Laplacian acting on co-exact $k$-forms on an $n=2k+1$-dimensional manifold are generically simple
if $k$ is odd and generically of multiplicity 2 if $k$ is even.


\subsection{Related work}\label{subsec:related1}

Bleeker and Wilson \cite{BW} studied eigenvalue multiplicity for the Laplace-Beltrami operator (the Hodge Laplacian on 0-forms) under conformal perturbations of the metric $g \rightarrow e^fg$, for $f \in C^\infty (M, \R)$ and proved generic simplicity of the eigenvalues.
More recently, Canzani \cite{canzani} studied the question of generic eigenvalue multiplicity for conformally covariant, elliptic self-adjoint  operators $P_g$ on
smooth sections of vector bundles over a compact Riemannian manifold $(M, g)$. Canzani proved that there is a residual set of functions in $C^\infty (M, \R)$ for which the corresponding operators $P_{e^fg}$ associated with the conformally deformed metrics $e^fg$ have simple nonzero eigenvalues.
The perturbation theory employed there, similar to that used in the present paper, depends crucially on the conformal covariance of the operators $P_g$. In related work, Jakobson and Strohmaier \cite{js1} studied quantum ergodicity for, among other operators, the Hodge Laplacian restricted to co-closed $k$-forms. In their study of quantum ergodicity for compact K\"ahler manifolds, Jacobson, Strohmaier, and Zelditch \cite[Remark 4.2]{jsz1} conjectured that the spectrum of the Hodge Laplacian restricted to primitive, co-closed $(p,q)$-forms is generically simple.

\subsection{Contents of the paper}

The Beltrami operator is studied in section 2. This is a skew-adjoint operator so the corresponding spectral problem is posed on the space of complex-valued 2-forms. It is shown in Theorem \ref{thm:beltrami-simple1} that its eigenvalues are generically simple. The relation between the eigenvalues of the Beltrami operator and Hodge Laplacian is discussed in section 3. The main result, Theorem \ref{thm:hodge-main1}, is proved in section 3, and states that the nonzero eigenvalues of the Hodge Laplacian acting on real-valued, co-exact 2-forms is generically two. In the last section, we discuss the general question of the generic multiplicity of the nonzero eigenvalues of the Hodge Laplacian acting on 2-forms over a 5-dimensional manifold.


\section{Generic simplicity of the eigenvalues of the Beltrami operator}\label{sec:beltrami-ev1}

The Beltrami operator $*_gd$ maps co-exact $2$-forms to co-exact $2$-forms on a $5$-dimensional manifold.
If $\omega$ is a co-exact $2$-form then it is easily found that
$$
\Delta_g^{(2)}\omega = \delta_g d\omega = -(*_gd)^2\omega.
$$
Furthermore, the Beltrami operator is skew-adjoint on the domain $H^1(M, \Lambda^2)$ in $L^2(M, \Lambda^2)$ with the inner product \eqref{eq:ip}.
Thus, in order to study the eigenvalues of the Beltrami operator, we consider the Beltrami operator on the space of complex-valued $2$-forms
$L^2 (M, \Lambda_{\C}^2)$.
Acting on its domain $H^1 ( M , \Lambda_\C^2)$, the Beltrami operator is skew-adjoint with purely imaginary eigenvalues.

We are interested in the multiplicities of the nonzero eigenvalues of the Beltrami operator restricted to the subspace of co-exact $2$-forms.
We define
 $$
 \mathcal{K}=\{u\in L^2(M,\Lambda^2)\,|\, du=0\},
 $$
 which is the set of all $L^2$ exact and harmonic 2-forms on $M$.  We will use $\perp_g$ to specify orthogonality with respect to the inner product \eqref{eq:ip}. By Hodge decomposition, $\mathcal{K}^{\perp_g}$ is the set of all $L^2$ co-exact 2-forms on $(M,g)$.  The spaces $\mathcal{K}$ and $\mathcal{K}^{\perp_g}$ consist of real 2-forms and will be used in section \ref{sec:hodge1}. In the present section in which we discuss the eigenvalue problem for the Beltrami operator, we will use the analogous spaces of complex-valued 2-forms, $\mathcal{K}_{\mathbb{C}}$ and $\mathcal{K}^{\perp_g}_{\mathbb{C}}$. The main result of this section is the generic simplicity of the eigenvalues of the Beltrami operator on co-exact $2$-forms.

\begin{thm}\label{thm:beltrami-simple1} The eigenvalues of the Beltrami operator $*_gd$ acting on the space $H^1(M,\Lambda_{\mathbb{C}}^2)\cap \mathcal{K}_{\mathbb{C}}^{\perp_g}$ are all simple for a residual set of $C^r$ metrics, for any $r \geq 2$.
\end{thm}

The proof of Theorem \ref{thm:beltrami-simple1} consists of a two parts. In the first,
we focus on one degenerate eigenvalue $i\lambda$ of $*_gd$. We prove that there is a real symmetric matrix $h$ so that the metric $g+ \epsilon h$ has a cluster of at least two nearby eigenvalues, converging to $i\lambda$ as $\epsilon \rightarrow 0$. Each will have multiplicity less than that of $i\lambda$.  In the second step, we prove that generically all eigenvalue multiplicities are removed using an inductive argument of Albert \cite[Theorems 1 and 2]{albert1} (see also Colin de Verdi\`ere, \cite[section 5]{cdv1}).


\subsection{Variation with respect to the metric}\label{subsec:variation1}

In this section, we compute the differential of the Beltrami operator $*_gd$ with respect to the metric $g$. Let $\mathcal{G}^r(M)$
denote the set of all $C^r$ metrics on the compact manifold $M$.
The space $\mathcal{S}^r(M)$ consists of all symmetric tensor fields of class $C^r$ and type $(0,2)$ and can be identified with the tangent space $T_g\mathcal{G}^r(M)$ at any $g\in\mathcal{G}^r(M)$.  Thus, $D(*d)_g(h)$ represents the variation of the Beltrami operator at the metric $g\in\mathcal{G}^r(M)$ in the direction of a $C^r$ symmetric $(0,2)$-tensor $h$. The trace of $h$ is given by $\tr_g h=g^{ij}h_{ij}$. The following lemma gives the local coordinate representation of $D(*d)_g(h)$ acting on an eigenform of the Beltrami operator.

\begin{lemma}\label{lemma:derivative}
Let $u\in H^1(M,\Lambda_{\mathbb{C}}^2)$ be an eigenform of $*_gd$ with eigenvalue $i\lambda$.  Then for any $h\in \mathcal{S}^r(M)$,	
\beq\label{eq:differential1}						
(D(*d)_g(h)u)_{ij} = i\lambda\left[-\frac{1}{2}(\tr_gh)u_{ij}+g^{mt}h_{ti}u_{mj}+g^{mt}h_{tj}u_{im}\right].
\eeq
\end{lemma}

\noindent
{\it Sketch of the proof.}
The proof of Lemma \ref{lemma:derivative} is computationally long. We provide an overview of the computations involved. Complete details are provided in \cite[Appendix A]{gier-thesis}.
First, we express the Beltrami operator in local coordinates by
$$
(*_gdu)_{ij} = \frac{1}{6} \varepsilon_{klmij}|g|^{1/2} g^{kn}g^{lp}g^{mq}\left(\frac{\partial u_{np}}{\partial x_q}-\frac{\partial u_{nq}}{\partial x_p}+\frac{\partial u_{pq}}{\partial x_n}\right).
$$
Next, using the formulas $$D(g^{ij})(h)=-h^{ij}\hspace{.2in}\mbox{and}\hspace{.2in}D(|g|^s)(h)=s|g|^s(\tr_g h) \hspace{.1in}\mbox{for } s>0,$$ we compute
  \begin{eqnarray*}
	(D(*d)_g(h)u)_{ij} &=& \frac{1}{6}\varepsilon_{klmij}|g|^{1/2}\left(\frac{\partial u_{np}}{\partial x_q}-\frac{\partial u_{nq}}{\partial x_p}+\frac{\partial u_{pq}}{\partial x_n}\right)\\
					   & & \times \left[\frac{1}{2}(\tr_gh)g^{kn}g^{lp}g^{mq}-g^{kn}g^{lp}h^{mq}-g^{kn}g^{mq}h^{lp}-g^{lp}g^{mq}h^{kn}\right].
   \end{eqnarray*}
Finally, we utilize the eigenvalue equation $*_gdu=i\lambda u$ to simplify the expression for $(D(*d)_g(h)u)_{ij}$. This results in the desired formula given in \eqref{eq:differential1}.\hfill $\Box$

\subsection{A density result}\label{subsec:density1}

The following density result states that any compactly-supported 2-form may be locally expressed in terms of a given non-vanishing form and a symmetric $(0,2)$-tensor.

\begin{lemma}\label{lem:density} Let $w\in C^r(M,\Lambda_{\mathbb{C}}^2)$, $r\geq 1$, and consider a compact subset $K\subset M\backslash w^{-1}(0)$.  Then for any
$v\in C^r(M,\Lambda_{\mathbb{C}}^2)$ with $\supp v\subset K$, there exists a symmetric complex $(0,2)$-tensor $t\in\mathcal{S}_{\mathbb{C}}^r(M)$ such that $v_{ij}=t_{ik}g^{kl}w_{lj}+w_{ik}g^{kl}t_{lj}$.
\end{lemma}

\noindent
{\it Sketch of the proof.} Let $w\in C^r(M,\Lambda_{\mathbb{C}}^2)$, let $K$ be a compact subset of $M\backslash w^{-1}(0)$, and let $v$ be any 2-form in $C^r(M,\Lambda_{\mathbb{C}}^2)$ with $\supp  v\subset K$.  To make the computations clearer, we will use matrix representations of the various forms and tensors.  The 2-forms $w$ and $v$ correspond to the antisymmetric $5\times 5$ matrices that we denote by $W$ and $V$, respectively. The 2-forms $g^{-1}$ and $t$ naturally correspond to the symmetric matrices denoted $G^{-1}$ and $T$.
The matrices $W,V,G^{-1},$ and $T$ are matrix-valued functions of $p\in M$.
The condition $v_{ij}=t_{ik}g^{kl}w_{lj}+w_{ik}g^{kl}t_{lj}$ for $1\leq i,j\leq 5$ translates into the matrix equation
$V=TG^{-1}W+WG^{-1}T$. Since $G^{-1}$ is a symmetric positive-definite matrix, it has a symmetric positive-definite square root $G^{-1/2}$.  We thus obtain the equivalent equation
\begin{eqnarray}\label{eqn:sylv}
\tilde{V} &=& \tilde{T}\tilde{W}+\tilde{W}\tilde{T},
\end{eqnarray}
where the matrices $\tilde{V}=G^{-1/2}VG^{-1/2}$ and $\tilde{W}=G^{-1/2}WG^{-1/2}$ are antisymmetric and $\tilde{T}=G^{-1/2}TG^{-1/2}$ is symmetric.

Let $\mathcal{M}$ denote the set of all $C^r$ $5\times 5$ matrix-valued functions on $M$. We define a linear operator
$L_{\tilde{W}}:\mathcal{M}\to\mathcal{M}$ by
\begin{eqnarray}\label{eqn:sylvester}
L_{\tilde{W}}(X)&:=& X\tilde{W}+\tilde{W}X.
\end{eqnarray}
Satisfying condition \eqref{eqn:sylv} amounts to finding a symmetric $\tilde{T}\in\mathcal{M}$ such that $L_{\tilde{W}}(\tilde{T})=\tilde{V}$.  The Sylvester equation
$L_{\tilde{W}}(X)=X\tilde{W}+\tilde{W}X=\tilde{V}$ has a unique solution if and only if
$\tilde{V}$ is orthogonal to $\ker L_{\tilde{W}}$ (see, for example \cite{Bhatia}). It is proved in \cite[Appendix C]{gier-thesis}
that each $E\in\ker L_{\tilde{W}}$ is symmetric.
By the antisymmetry of $\tilde{V}$, the matrix inner product of $\tilde{V}$ with each $E\in \ker L_{\tilde{W}}$ is
		\begin{eqnarray*}
									 E\cdot \tilde{V} &=& \sum_{i,j=1}^5 \overline{e_{ij}}\tilde{v}_{ij} \\
																		&=& \sum_{i<j} \overline{e_{ij}}\tilde{v}_{ij} + \sum_{i>j} \overline{e_{ij}}\tilde{v}_{ij} \\
				                            &=& \sum_{i<j} \overline{e_{ij}}\tilde{v}_{ij} + \sum_{i>j} \overline{e_{ji}}(-\tilde{v}_{ji}) \\
				                            &=& \sum_{i<j} \overline{e_{ij}}\tilde{v}_{ij} - \sum_{i<j} \overline{e_{ij}}\tilde{v}_{ij} \hspace{.2in}\mbox{(reindexing)}\\
				                            &=& 0.
		\end{eqnarray*}
Since $\tilde{V}$ is orthogonal to $\ker L_{\tilde{W}}$,
there exists an $X\in\mathcal{M}$ such that $\tilde{V}=X\tilde{W}+\tilde{W}X$ on $K$.  From the antisymmetry of $\tilde{V}$ and $\tilde{W}$, one easily shows that $X^T\tilde{W}+\tilde{W}X^T = \tilde{V}$,
so that $X^T$ solves the same equation as $X$.  Thus, we define $\tilde{T}$ to be the symmetrization $\tilde{T}=\frac{1}{2}[X+X^T]$. Hence, $T=G^{1/2}\tilde{T} G^{1/2}$ is a symmetric $C^r$ matrix-valued function such that $V=TG^{-1}W+WG^{-1}T$. We thus obtain from $T$ the desired symmetric complex $(0,2)$-tensor $t\in\mathcal{S}_{\mathbb{C}}^r(M)$. \hfill $\Box$


\subsection{Eigenvalue perturbation theory}\label{subsec:eigenvalue-pert1}

To establish the generic simplicity of the eigenvalues of the Beltrami operator, we use standard results from perturbation theory as discussed in Rellich \cite[chapter II, section 5, Theorem 3]{Rellich} and Kato \cite{Kato}.  In particular, observe that the skew-adjointness of the Beltrami operator $*_gd$ when $n=5$ and $k=2$ implies that the operator $i\hspace{-.03in}*_gd:H^1(M,\Lambda_{\mathbb{C}}^2)\cap\mathcal{K}_{\mathbb{C}}^{\perp_g}\to \mathcal{K}_{\mathbb{C}}^{\perp_g}$ is self-adjoint with respect to the metric $g$ and has real, isolated eigenvalues of finite multiplicity. We consider perturbations of the metric $g \rightarrow g(\epsilon) := g + \epsilon h$ so the norm, and hence the Hilbert space, depends on $\epsilon$. We map these spaces to the $\epsilon$-independent Hilbert space $L^2 (M, \Lambda^2_{\C})$. We define a unitary operator $U_\epsilon :  L^2 (M, \Lambda^2_{\C}) \rightarrow L^2 (M, \Lambda^2_{\C}, g(\epsilon))$ by
$$
U_\epsilon \omega = \left( \frac{ \det g }{\det g(\epsilon) } \right)^{1/4} \omega,
$$
for any two-form $\omega \in L^2 (M, \Lambda^2_{\C})$. Then the Beltrami operator $\mathcal{D}_\epsilon := U_\epsilon^{-1} ( *_{g(\epsilon)} d) U_\epsilon$
acts on $L^2 (M, \Lambda^2_{\C})$ and is unitarily equivalent to the Beltrami operator $*_{g (\epsilon)} d$.
Note that $\mathcal{D}_0 = *_g d$. Furthermore, the set of co-exact two-forms $\mathcal{K}^{\perp_{g(\epsilon)}}$
in $L^2 (M, \Lambda^2_{\C}, g(\epsilon))$ maps to the $\mathcal{D}_\epsilon$-invariant
subspace $\tilde{\mathcal{K}} ^{\perp_{g(\epsilon)}} \subset L^2 (M,  \Lambda^2_{\C}, g)$.

In this setting, we have the following perturbation theorem for linear perturbations of the metric.

\begin{thm}\label{thm:perttheoremh}
Let $\lambda$ be an eigenvalue of $i\hspace{-.03in}*_gd:H^1(M,\Lambda_{\mathbb{C}}^2)\cap\mathcal{K}_{\mathbb{C}}^{\perp_g}\to \mathcal{K}_{\mathbb{C}}^{\perp_g}$ of multiplicity $m$, and let $g(\epsilon)=g+\epsilon h$ for some $h\in S^r(M)$.  Then there are $m$ functions $\ell^h_1(\epsilon),\ldots,\ell^h_m(\epsilon)$ real-analytic at $\epsilon=0$, and $m$ functions $U^h_1(\epsilon),\ldots,U^h_m(\epsilon) \in
L^2 (M, \Lambda_\C^2)$, analytic in $H^1(M,\Lambda_{\mathbb{C}}^2)$ at $\epsilon=0$, such that the following conditions hold:
\begin{enumerate}
  \item $\ell^h_j(0)=\lambda$ for $j=1,\ldots,m$;
    \item $i \mathcal{D}_\epsilon U^h_j(\epsilon) = \ell^h_j(\epsilon)U^h_j(\epsilon)$ for $j=1,\ldots,m$;
    \item For $\epsilon$ in a small enough neighborhood of $0$, $\{U^h_1(\epsilon),\ldots, U^h_m(\epsilon)\}$ is an orthonormal set in $H^1(M,\Lambda_{\mathbb{C}}^2)\cap \tilde{\mathcal{K}}_{\mathbb{C}}^{\perp_{g(\epsilon)}}$;
    \item For every open interval $(a,b)\subset\mathbb{R}$ such that $\lambda$ is the only eigenvalue of $i\hspace{-.03in}*_gd$ in $[a,b]$, there are exactly $m$ eigenvalues (counting multiplicity) $\ell^h_1(\epsilon),\ldots,\ell^h_m(\epsilon)$ of $i\hspace{-.03in}*_{g(\epsilon)}d$ in $(a,b)$, for $\epsilon$ sufficiently small.
\end{enumerate}
\end{thm}

It will be convenient for the calculation in section \ref{subsec:beltrami-simple1} to write the eigenvalue equation
in the second point of Theorem \ref{thm:perttheoremh} in the following form. Since
\beq\label{eq:pert1}
i \mathcal{D}_\epsilon U^h_j(\epsilon) = i U_\epsilon^{-1} (\hspace{-.03in}*_{g(\epsilon)}d ) ( U_\epsilon U^h_j(\epsilon)),
\eeq
if we let $\tilde{U}_j^h (\epsilon) := U_\epsilon U^j_h (\epsilon)$, we have
\beq\label{eq:pert2}
i \hspace{-.03in}*_{g(\epsilon)}d \tilde{U}^h_j(\epsilon) = \ell^h_j(\epsilon) \tilde{U}^h_j(\epsilon),
\eeq
for $j=1,\ldots,m$. These eigenforms $\tilde{U}_j^h (\epsilon)$ belong to $L^2 (M, \Lambda_\C^2, g(\epsilon))$.


\subsection{Proof of Theorem \ref{thm:beltrami-simple1} for the Beltrami operator}\label{subsec:beltrami-simple1}

We combine the perturbation result with the topological arguments of Albert to prove Theorem \ref{thm:beltrami-simple1}.
\newline
\noindent
\textit{Proof.}\\
\noindent
1. The setting. For a metric $g\in\mathcal{G}^r(M)$, we label the eigenvalues $i\lambda_n (g)$ of the Beltrami operator $*_gd$ so that
$$\lambda_{n+1}^2(g)\geq \lambda_n^2(g).$$
We define the following subsets of $\mathcal{G}^r(M)$, the metrics on $M$:
$$\Gamma_\infty := \{g\in\mathcal{G}^r(M)\,|\,\mbox{ all eigenvalues of }*_gd|_{H^1(M,\Lambda_{\mathbb{C}}^2)\cap \mathcal{K}_{\mathbb{C}}^{\perp_g}}\mbox{ are simple}\}$$
and
$$\Gamma_n :=\{g\in\mathcal{G}^r(M)\,|\,\mbox{ the first }n\mbox{ eigenvalues of }*_gd|_{H^1(M,\Lambda_{\mathbb{C}}^2)\cap \mathcal{K}_{\mathbb{C}}^{\perp_g}}\mbox{ are simple}\}.
$$
These subsets are nested so that
$$
\Gamma_\infty \subset\cdots\subset \Gamma_n\subset\Gamma_{n+1}\subset\cdots\subset\Gamma_1\subset\Gamma_0=\mathcal{G}^r(M),
$$
and $$\Gamma_\infty =\bigcap_{n=0}^\infty \Gamma_n.$$
By the stability of simple eigenvalues under small perturbations of the metric, each set $\Gamma_n$ is open in $\mathcal{G}^r(M)$.  Thus, to prove that $\Gamma_\infty$ is residual in $\mathcal{G}^r(M)$, it is sufficient to show that $\Gamma_{n+1}$ is dense in $\Gamma_{n}$ for all $n=0,1,2,\ldots$.

\noindent
2. The density argument. Let $g\in \Gamma_n$ so that the first $n$ eigenvalues of $$*_gd:H^1(M,\Lambda_{\mathbb{C}}^2)\cap\mathcal{K}_{\mathbb{C}}^{\perp_g}\to\mathcal{K}_{\mathbb{C}}^{\perp_g}$$ are simple.  Suppose that the $(n+1)^{\rm st}$ eigenvalue $i\lambda\neq 0$ of $*_gd$ has multiplicity $m$, and define $g(\epsilon)=g+\epsilon h$ for some $h\in S^r(M)$. Theorem \ref{thm:perttheoremh} implies there are $m$ functions $\ell^h_1(\epsilon),\ldots,\ell^h_m(\epsilon)$ real-analytic at $\epsilon=0$, and $m$ functions $U^h_1(\epsilon),\ldots,U^h_m(\epsilon)$ analytic in $H^1(M,\Lambda_{\mathbb{C}}^2)$ at $\epsilon=0$ such that the conditions of Theorem \ref{thm:perttheoremh} hold.
When $\epsilon=0$, each set $\{U^h_1(0),\ldots, U^h_m(0)\}$ forms an orthonormal basis of the eigenspace $E(*_gd,i\lambda)$. This basis may depend on the choice of $h\in\mathcal{S}^r(M)$ in the linear perturbation of the metric $g(\epsilon)=g+\epsilon h$.

\noindent
3. Variation with respect to the metric.
We differentiate the eigenvalue equation
\begin{eqnarray*}
\displaystyle *_{g(\epsilon)}d {\tilde U}^h_j(\epsilon)&=& i\ell^h_j(\epsilon) {\tilde U}^h_j(\epsilon),
\end{eqnarray*}
where ${\tilde U}^h_j(\epsilon) = U_\epsilon {U}^h_j(\epsilon) \in L^2(M,\Lambda_{\mathbb{C}}^2, g(\epsilon)$,
with respect to $\epsilon$ and evaluate at $\epsilon=0$ to obtain
\beq\label{eq:diff-ev-eqn1}
D(*d)_g(h)U^h_j(0)+*_gd({\tilde U}^h_j)'(0) = i(\ell^h_j)'(0)U^h_j(0)+i\ell^h_j(0)({\tilde U}^h_j)'(0),
\eeq
where $({\tilde U}^h_j)'(0) \in L^2 (M,\Lambda_{\mathbb{C}}^2)$ due to the analyticity in Theorem \ref{thm:perttheoremh}.
Introducing the notation $u^h_j=U^h_j(0)$, we simplify \eqref{eq:diff-ev-eqn1} to
\beq\label{eqn:perderh}
D(*d)_g(h)u^h_j+(*_gd-i\lambda)({\tilde U}^h_j)'(0) = i(\ell^h_j)'(0)u^h_j.
\eeq
Since $\{u_1^h,\ldots,u_m^h\}$ is an orthonormal basis of the eigenspace
$E(*_gd,i\lambda)$, we take the inner product of \eqref{eqn:perderh} with another eigenform $u^h_k$. This results in
\beq
i(\ell^h_j)'(0)(u^h_j,u^h_k)_g = (D(*d)_g(h)u^h_j,u^h_k)_g+((*_gd-i\lambda)({\tilde U}^h_j)'(0),u^h_k)_g .
\eeq
The last term on the right vanishes due to the skew-adjointness of the Beltrami operator and the eigenvalue equation.
Consequently, we obtain
\beq
i(\ell^h_j)'(0)\delta_{jk} = (D(*d)_g(h)u^h_j,u^h_k)_g.
\eeq
We may express the inner product $(D(*d)_g(h)u^h_j,u^h_k)_g$ in local coordinates using Lemma \ref{lemma:derivative} to obtain
\begin{eqnarray}\label{eq:variation1}
(\ell^h_j)'(0)\delta_{jk}  &=& \frac{\lambda}{2}\int g^{pr}g^{qs}\left[-\frac{1}{2}(\tr_gh)(u^h_j)_{pq}+g^{lt}h_{tp}(u^h_j)_{lq}+g^{lt}h_{tq}(u^h_j)_{pl}\right]
                                                            \overline{(u^h_k)_{rs}}\,d\mu_g. \nonumber \\
                                                             & &
\end{eqnarray}
We define a bilinear form $S:\mathcal{S}^r_{\mathbb{C}}(M)\times L^2(M,\Lambda_{\mathbb{C}}^2)\to L^2(M,\Lambda_{\mathbb{C}}^2)$ by
  \beq\label{eq:bilinear1}
  [S(h, w)]_{pq}=-\frac{1}{2}(\tr_gh)w_{pq}+g^{lt}h_{tp}w_{lq}+g^{lt}h_{tq}w_{pl}
  \eeq
    for $h\in\mathcal{S}^r_{\mathbb{C}}(M)$ and $w\in  L^2(M,\Lambda_{\mathbb{C}}^2)$.
We may then express \eqref{eq:variation1} more concisely as
    \begin{eqnarray}\label{eqn:orthwh}
        (\ell^h_j)'(0)\delta_{jk} &=& \lambda(S(h, u^h_j),u^h_k)_g.
    \end{eqnarray}

\noindent
4. Change of basis. Our goal is to show that there exists an
$h\in\mathcal{S}^r(M)$ such that $$(\ell^h_j)'(0)\neq (\ell^h_k)'(0)$$ for some pair $j,k\in\{1,\ldots,m\}$.
This fact implies that under the metric perturbation $g(\epsilon)=g+\epsilon h$ for $\epsilon$
sufficiently small, the perturbed eigenvalues $i\ell^h_j(\epsilon)$ and $i\ell^h_k(\epsilon)$ of $*_{g(\epsilon)}d$ are distinct.
While $i\ell^h_j(\epsilon)$ and $i\ell^h_k(\epsilon)$ are not guaranteed to be simple, they each have multiplicity less than $m$.
 To this end, assume to the contrary that $(\ell^h_j)'(0) = (\ell^h_k)'(0)$ for all $h\in\mathcal{S}^r(M)$ and all $j,k\in\{1,\ldots,m\}$.  By \eqref{eqn:orthwh}, this assumption implies
\begin{eqnarray}
(S(h, u^h_j),u^h_j)_g &=& (S(h, u^h_k),u^h_k),\hspace{.1in} 1\leq j,k\leq m \label{eqn:scond1}\\
(S(h, u^h_j),u^h_k)_g &=& 0,\hspace{.1in} j\neq k \label{eqn:scond2}
\end{eqnarray}
for all $h\in \mathcal{S}^r(M)$.  As previously noted, each set $\{u_1^h,\ldots,u_m^h\}$ forms an orthonormal basis of $E(*_gd, i\lambda)$, but we cannot assume that $u^{h_1}_j=u^{h_2}_j$ when $h_1\neq h_2$.  Let us therefore fix an orthonormal basis $\{u_1,\ldots,u_m\}$ of $E(*_gd,i\lambda)$.  For a given $h\in \mathcal{S}^r(M)$, we write each $u_j$ in terms
of the basis elements $\{u_1^h,\ldots,u_m^h\}$ as $u_j = \sum_{\ell = 1}^m c_{j, \ell} u_\ell^h$,
for constants $c_{j,\ell} \in \mathbb{C}$.  The fact that $\{u_1,\ldots,u_m\}$ and $\{u_1^h,\ldots,u_m^h\}$ are both orthonormal bases of $E(*_gd,i\lambda)$ implies
\beq\label{eqn:cd}  
\delta_{jk} = (u_j,u_k)_g = \sum_{\ell = 1}^m c_{j, \ell} \overline{c_{k,\ell}}.
\eeq
Combining \eqref{eqn:cd} with \eqref{eqn:scond1} and \eqref{eqn:scond2} yields
\begin{eqnarray*}
(S(h, u_j),u_k)_g &=& c_{j,1}(S(h, u^h_1),u_k)_g+\cdots+c_{j,m}(S(h, u^h_m),u_k)_g \\
                  &=& c_{j,1}\overline{c_{k,1}}(S(h, u^h_1),u^h_1)_g+\cdots+c_{j,m}\overline{c_{k,m}}(S(h, u^h_m),u^h_m)_g\\
                  &=& (c_{j,1}\overline{c_{k,1}}+\cdots+c_{j,m}\overline{c_{k,m}})(S(h, u^h_j),u^h_j)_g\\
                  &=& \delta_{jk}(S(h, u^h_j),u^h_j)_g.
\end{eqnarray*}
Thus, for all $h\in\mathcal{S}^r(M)$, the elements in the orthonormal basis $\{u_1,\ldots,u_m\}$ satisfy
\bea\label{eq:basis-rel1}
(S(h, u_j),u_j)_g &=& (S(h, u_k),u_k)_g, \hspace{.1in} 1\leq j,k\leq m \nonumber \\
(S(h, u_j),u_k)_g &=& 0, \hspace{.1in}j\neq k.
\eea

\noindent
5. Extension to $\mathcal{S}^r_{\mathbb{C}}(M)$.  For any $T\in \mathcal{S}^r(M)$, we define $h_T$ by
    \beq\label{eqn:htwh}
        h_T = T-(\tr_g T)g.
    \eeq
With this choice, the bilinear form $S$ defined in \eqref{eq:bilinear1} becomes
	\begin{eqnarray*}
				[S(h_T, u_j)]_{pq}&=& -\frac{1}{2}[(\tr_g T)-5(\tr_g T)](u_j)_{pq}+g^{lt}[T_{tp}-(\tr_g T)g_{tp}](u_j)_{lq}\\
                             & & +g^{lt}[T_{tq}-(\tr_g T)g_{tq}](u_j)_{pl} \\
							 &=& 2(\tr_g T)(u_j)_{pq}+g^{lt}T_{tp}(u_j)_{lq}-(\tr_g T)(u_j)_{pq}+g^{lt}T_{tq}(u_j)_{pl}\\
                             & & -(\tr_g T)(u_j)_{pq}\\
							 &=& T_{pt}g^{tl}(u_j)_{lq}+(u_j)_{pl}g^{lt}T_{tq}.
		\end{eqnarray*}
By decomposing a complex symmetric $(0,2)$-tensor $T\in \mathcal{S}^r_{\mathbb{C}}(M)$ into $T=T_1+iT_2$ for $T_1,T_2\in \mathcal{S}^r(M)$, the linearity of $h_T$ in $T$ \eqref{eqn:htwh} and relations \eqref{eq:basis-rel1} for real $T$
imply $(S(h_{T},u_j), u_j)_g = (S(h_T, u_k),u_k)_g$,
for all $T\in \mathcal{S}_{\mathbb{C}}^r(M)$.  Likewise, we obtain
\begin{eqnarray}\label{eqn:what2}
(S(h_T,u_j),u_k)_g = 0,&& j\neq k
\end{eqnarray}
for all complex tensors $T\in S_{\mathbb{C}}^r(M)$.

\noindent
6. Unique continuation principle.  Without loss of generality, we fix $j=1$ and $k=2$. Equation \eqref{eqn:what2} implies
\begin{eqnarray}\label{eqn:quicker}
(S(h_T,u_1),u_2)_g &=& 0
\end{eqnarray}
for all $T\in \mathcal{S}_{\mathbb{C}}^r(M)$.  We apply Lemma \ref{lem:density} with $w = u_1$.
It follows from $(*_gd-i\lambda)u_1=0$ and the co-exactness of $u_1$ that
$\Delta_g^{(2)}u_1=-(*_gd)^2u_1=\lambda^2 u_1$,
 so that $u_1$ is an eigenform of the Hodge Laplacian $\Delta^{(2)}_g$ with eigenvalue $\lambda^2$.  The unique continuation principle
then states that $u_1$ cannot vanish in any open subset of $M$ \cite{Aronszajn,uniquecont}. Consequently, the set
$$\mathscr{S}=\{S(h_{T},u_1)\,|\, T\in\mathcal{S}_{\mathbb{C}}^r(M)\}$$
is dense in $L^2(M,\Lambda_{\mathbb{C}}^2)$ by Lemma \ref{lem:density}.  Since \eqref{eqn:quicker} implies $u_2$ is orthogonal to the dense set $\mathscr{S}$, we obtain $u_2=0$ on $M$, contradicting the fact that $u_2$ is a normalized eigenform.

\noindent
7. Conclusion of the proof. By the above, there exists an $h\in \mathcal{S}^r(M)$ such that
$(\ell^h_j)'(0) \neq (\ell^h_k)'(0)$ for some $j,k\in\{1,\ldots,m\}$. Consequently, for all $\epsilon > 0$ small, the $n+1^{\rm st}$ eigenvalue of the Beltrami operator $*_{g(\epsilon)} d$ has multiplicity at most $m-1$, and the first $n$ eigenvalues remain simple.   Repeating the above argument as necessary, we obtain a metric $g(\epsilon)=g+\epsilon h$ in $\Gamma^{n+1}$ for $\epsilon$ sufficiently small.  Since $g(\epsilon)$ can be taken arbitrarily close to $g$ in the $C^r$ topology, we conclude that $\Gamma^{n+1}$ is dense in $\Gamma^n$. Additionally, each $\Gamma^n$ is open in $\mathcal{G}^r(M)$, so we infer that
$$\Gamma_\infty =\bigcap_{n=1}^\infty \Gamma_n$$
is residual $\mathcal{G}^r(M)$.  Thus, for a residual set of metrics $\Gamma_\infty \subset\mathcal{G}^r(M)$, the Beltrami operator acting on $H^1(M,\Lambda_{\mathbb{C}}^2)\cap \mathcal{K}_{\mathbb{C}}^{\perp_g}$ has only simple eigenvalues.\hfill $\Box$


\section{The Hodge Laplacian on co-exact $2$-forms}\label{sec:hodge1}

In this section, we apply the results on the generic multiplicities of the Beltrami operator to the study of the eigenvalue multiplicities of the Hodge Laplacian acting on real co-exact 2-forms.


\subsection{Relation to the eigenvalues of the Beltrami operator}\label{subsec:relation1}

In order to determine the generic eigenvalue multiplicities of the Hodge Laplacian on co-exact 2-forms on a 5-manifold, we must determine the relationship between the eigenvalues and eigenforms of the Hodge Laplacian and those of the Beltrami operator.  Our next two lemmas hold in the more general setting of $n=4\ell+1$ and $k=2\ell$ for some $\ell\in\mathbb{N}$ and in particular apply when $n=5$ and $k=2$.

\begin{lemma}\label{lem:bab}
Let $M$ be a manifold of dimension $n=4\ell+1$ for some $\ell\in\mathbb{N}$, and let $k=2\ell$.  Let $\omega=\alpha+i\beta$ be a nonzero complex $k$-form with $\alpha,\beta\in H^1(M,\Lambda^k)$.  Then $*_gd\omega=i\lambda\omega$ if and only if
\begin{eqnarray}\label{eqn:ab}
	*_gd\alpha &=& -\lambda\beta \hspace{.2in}\mbox{and}\hspace{.2in} *_gd\beta=\lambda\alpha.
\end{eqnarray}
\end{lemma}

\textit{Proof.} First, suppose that $\omega=\alpha+i\beta\in H^1(M,\Lambda_{\mathbb{C}}^k)$ solves $*_gd\omega=i\lambda\omega$.  Then
\begin{eqnarray*}
	           *_gd\omega &=& i\lambda\omega \\
	  *_gd(\alpha+i\beta) &=& i\lambda(\alpha+i\beta)\\
	*_gd\alpha+i*_gd\beta &=& -\lambda\beta+i\lambda\alpha,
\end{eqnarray*}
so equating real and imaginary parts yields \eqref{eqn:ab}.

Conversely, suppose that $\alpha,\beta\in H^1(M,\Lambda^k)$ satisfy \eqref{eqn:ab}, and let $\omega=\alpha+i\beta$.  Then
$$*_gd\omega=*_gd\alpha+i*_gd\beta=-\lambda\beta+i\lambda\alpha=i\lambda(\alpha+i\beta)=i\lambda\omega$$
so that $\omega$ is an eigenfunction of $*_gd$ with eigenvalue $i\lambda$. \hfill $\Box$

\begin{remark}\label{n:independent}
 It is important to recognize that condition \eqref{eqn:ab} implies that $\alpha$ and $\beta$ are nonzero, linearly independent forms over $\mathbb{R}$.  To see this, observe that $\beta=c\alpha$ implies
$$\beta=c\alpha=\frac{c}{\lambda}*_gd\beta=\frac{c^2}{\lambda}*_gd\alpha=-c^2\beta,$$ which gives $c=\pm i$ in contradiction to $c\in\mathbb{R}$. Even more notably,
$$(\alpha,\beta)_g = \frac{1}{\lambda}(*_gd\beta,\beta)_g = -\frac{1}{\lambda}(\beta,*_gd\beta)_g=-(\beta,\alpha)_g$$ reveals that
\beq\label{eqn:aborth}
(\alpha,\beta)_g = 0.
\eeq
\end{remark}

The next lemma follows from our observations in Lemma \ref{lem:bab}.

\begin{lemma}\label{lem:bhl}
Let $M$ be a manifold of dimension $n=4\ell+1$ for some $\ell\in\mathbb{N}$, and let $k=2\ell$.  Let $\alpha,\beta\in H^2(M,\Lambda^k)\cap \mathcal{K}^{\perp_g}$.  If $\omega=\alpha+i\beta$ is an eigenform of the Beltrami operator $*_gd$ with eigenvalue $i\lambda$, then both $\alpha$ and $\beta$ are eigenforms of the Hodge Laplacian $\Delta_g^{(k)}$ with eigenvalue $\lambda^2$.
\end{lemma}

\textit{Proof.} Let $\alpha,\beta\in H^2(M,\Lambda^k)\cap \mathcal{K}^{\perp_g}$, and suppose $\omega=\alpha+i\beta$ is an eigenform of $*_gd$ with eigenvalue $i\lambda$.  By Lemma \ref{lem:bab}, $\alpha$ and $\beta$ satisfy
\begin{eqnarray*}
	*_gd\alpha &=& -\lambda\beta \hspace{.2in}\mbox{and}\hspace{.2in} *_gd\beta=\lambda\alpha.
\end{eqnarray*}
Since $\alpha$ is a co-exact form, $n=4\ell+1$ is odd, and $k=2\ell$ is even,
$$\Delta_g^{(k)}\alpha=-(*_gd)^2\alpha=\lambda*_gd\beta=\lambda^2\alpha.$$
Similarly,
$$\Delta_g^{(k)}\beta=-(*_gd)^2\beta=-\lambda*_gd\alpha=\lambda^2\beta$$
so that $\alpha$ and $\beta$ are both eigenforms of $\Delta_g^{(k)}$ with eigenvalue $\lambda^2$. \hfill $\Box$

\subsection{Proof of Theorem \ref{thm:hodge-main1}}

\textit{Proof.} By Theorem \ref{thm:beltrami-simple1}, there exists a residual set $\Gamma$ of $C^r$ metrics on $M$ such that the eigenvalues of the Beltrami operator $*_gd$ acting on $H^1(M,\Lambda_{\mathbb{C}}^2)\cap \mathcal{K}_{\mathbb{C}}^{\perp_g}$ are all simple.  Take $g\in \Gamma$, and consider an eigenvalue $\lambda^2>0$ of the restriction of $\Delta_g^{(2)}$ to co-exact 2-forms.  Let $\eta\in H^2(M,\Lambda^2)\cap \mathcal{K}^{\perp_g}$ be an eigenform of $\Delta_g^{(2)}$ with eigenvalue $\lambda^2$ so that
\begin{eqnarray}\label{eqn:eb}
-(*_gd)^2\eta &=& \lambda^2\eta.
\end{eqnarray}
Now, since $*_gd$ maps $H^2(M,\Lambda^2)\cap \mathcal{K}^{\perp_g}$ to $H^1(M,\Lambda^2)\cap \mathcal{K}^{\perp_g}$, we have $*_gd\eta=\lambda\zeta$ for some $\zeta\in H^1(M,\Lambda^2)\cap\mathcal{K}^{\perp_g}$.  Equation \eqref{eqn:eb} then yields
\begin{eqnarray*}
-*_gd(\lambda\zeta) &=& \lambda^2\eta\\
          *_gd\zeta &=& -\lambda\eta
\end{eqnarray*}
so that $\zeta$ is in fact contained in $H^2(M,\Lambda^2)\cap \mathcal{K}^{\perp_g}$.  Since $\eta$ and $\zeta$ together satisfy condition \eqref{eqn:ab}, Lemma \ref{lem:bab} implies that $\zeta+i\eta$ is an eigenform of $*_gd$ with eigenvalue $i\lambda$.  It follows from Lemma \ref{lem:bhl} that $\zeta$ is also an eigenform of $\Delta_g^{(2)}$ with eigenvalue $\lambda^2$.  As mentioned in \ref{n:independent}, the eigenforms $\eta$ and $\zeta$ are linearly independent, indicating that the eigenvalue $\lambda^2$ of $\Delta_g^{(2)}$ has multiplicity of at least 2.

To prove that $\lambda^2$ has a multiplicity of precisely 2, suppose that $\Delta_g^{(2)}\tau=\lambda^2\tau$ for some $\tau\in H^2(M,\Lambda^2)\cap \mathcal{K}^{\perp_g}$.  By our previous argument, there must exist a co-exact 2-form $\xi\in H^2(M,\Lambda^2)\cap \mathcal{K}^{\perp_g}$ such that $\xi+i\tau$ is an eigenform of $*_gd$ with eigenvalue $i\lambda$.  Since $g$ is contained in the residual set $\Gamma$, the eigenvalue $i\lambda$ is simple.  Thus, $\xi+i\tau$ must be a complex multiple of the eigenform $\zeta+i\eta$; that is,
\begin{eqnarray}\label{eqn:xz}
\xi+i\tau\hspace{.08in}=\hspace{.08in}(a+ib)(\zeta+i\eta)&=&(a\zeta-b\eta)+i(b\zeta+a\eta)
\end{eqnarray} for some $a+ib\in\mathbb{C}$.  Equating the imaginary parts of equation \eqref{eqn:xz} gives
$$\tau=b\zeta+a\eta$$
so that $\tau$ is a linear combination of the eigenforms $\eta$ and $\zeta$ of $\Delta_g^{(2)}$.  Thus, $\lambda^2$ has multiplicity 2.  We therefore conclude that for a residual set of metrics $\Gamma\subset \mathcal{G}^r(M)$, all eigenvalues of the restriction of the Hodge Laplacian $\Delta_g^{(2)}$ to $H^2(M,\Lambda^2)\cap \mathcal{K}^{\perp_g}$ have multiplicity 2. \hfill $\Box$

\vspace{.2in}

As an immediate consequence of the commutativity of the Hodge Laplacian and the exterior differential operator, we obtain an analogous result for exact 3-forms.

\begin{follow}\label{cor:ethree}
Let $M$ be a closed 5-manifold, and let $r$ be an integer, $r\geq 2$.  There exists a residual subset $\Gamma$ of the space of all $C^r$ metrics on $M$ such that, for all $g\in\Gamma$, the eigenvalues of the restriction of the Hodge Laplacian $\Delta_g^{(3)}$ to the space of exact forms in $H^2(M,\Lambda^3)$ have multiplicity 2.
\end{follow}


\section{Conclusions: Multiplicities of Hodge eigenvalues on $5$-manifolds}\label{sec:conclusions1}

Let us summarize the situation concerning the generic multiciplities of the nonzero eigenvalues of the Hodge Laplacian $\Delta_g^{(k)}$ on a closed 5-manifold. Uhlenbeck's reults \cite{uhlenbeck1}
ensures the generic simplicity of the nonzero eigenvalues of following operators:
\begin{enumerate}
	\item[(i)] $\Delta_g^{(0)}$;
	\item[(ii)] $\Delta_g^{(1)}$ restricted to exact 1-forms;
	\item[(iii)] $\Delta_g^{(4)}$ restricted to co-exact 4-forms;\
	\item[(iv)] $\Delta_g^{(5)}$.
\end{enumerate}
The first statement, the generic simplicity of the nonzero eigenvalues of $\Delta_g^{(0)}$, is Theorem 8 of \cite{uhlenbeck1}. The second statement
follows from this and the fact that an exact 1-form has the form $df$, for a function $f$. The last two statements follow from Hodge duality.
Moreover, Theorem \ref{thm:hodge-main1} and Corollary \ref{cor:ethree} of the present work assert
that there exists a residual set of $C^r$ metrics such that the operators
\begin{enumerate}
	\item[(v)] $\Delta_g^{(2)}$ restricted to co-exact 2-forms,
	\item[(vi)] $\Delta_g^{(3)}$ restricted to exact 3-forms
\end{enumerate}
have eigenvalues of multiplicity 2.
In order to completely characterize the generic nonzero eigenvalue multiplicities of the Hodge Laplacian on a closed 5-manifold, we also need information about the eigenspaces of the operators
\begin{enumerate}
	\item[(vii)] $\Delta_g^{(1)}$ restricted to co-exact 1-forms,
	\item[(viii)] $\Delta_g^{(2)}$ restricted to exact 2-forms,
	\item[(ix)] $\Delta_g^{(3)}$ restricted to co-exact 3-forms,
	\item[(x)] $\Delta_g^{(4)}$ restricted to exact 4-forms.
\end{enumerate}
Since operators (vii)-(x) have isomorphic eigenspaces, it suffices to determine the eigenvalue multiplicities of the Hodge Laplacian restricted to co-exact 1-forms.  It is unclear how to best approach this problem.  On a 5-manifold, the Beltrami operator only maps the space of 2-forms to itself and hence only has eigenvalues when acting on 2-forms.  Thus, the eigenvalue multiplicities of the Beltrami operator will not give insight into the eigenvalues of $\Delta_g^{(1)}$ on co-exact 1-forms. A direct perturbation approach is possible but to obtain results similar to those in section \ref{sec:beltrami-ev1}, and Lemma \ref{lem:density} in particular, for the Hodge Laplacian, would be calculationally intensive.


\end{document}